**RESEARCH**     **Open Access**

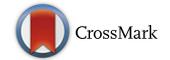

# On a unified framework for linear nuisance parameters

Yongchang Hu* and Geert Leus

**Abstract**

Estimation problems in the presence of deterministic linear nuisance parameters arise in a variety of fields. To cope with those, three common methods are widely considered: (1) jointly estimating the parameters of interest and the nuisance parameters; (2) projecting out the nuisance parameters; (3) selecting a reference and then taking differences between the reference and the observations, which we will refer to as "differential signal processing." A lot of literature has been devoted to these methods, yet all follow separate paths.

Based on a unified framework, we analytically explore the relations between these three methods, where we particularly focus on the third one and introduce a general differential approach to cope with multiple distinct nuisance parameters. After a proper whitening procedure, the corresponding best linear unbiased estimators (BLUEs) are shown to be all equivalent to each other. Accordingly, we unveil some surprising facts, which are in contrast to what is commonly considered in literature, e.g., the reference choice is actually not important for the differencing process. Since this paper formulates the problem in a general manner, one may specialize our conclusions to any particular application. Some localization examples are also presented in this paper to verify our conclusions.

## 1 Introduction

The problem of estimating unknown parameters of interest $\mathbf{x} \in \mathbb{R}^{L \times 1}$ observed through a linear transformation $\mathbf{H} \in \mathbb{R}^{N \times L}$ ($N > L$), and corrupted by additive noise $\mathbf{n} \in \mathbb{R}^{N \times 1}$, has been well studied and considered in a wide variety of fields [1]. However, the observations $\mathbf{y} \in \mathbb{R}^{N \times 1}$ are sometimes also influenced by unknown linear nuisance parameters, denoted by $\mathbf{u} \in \mathbb{R}^{M \times 1}$ which enter $\mathbf{y}$ through the linear transformation $\mathbf{G} \in \mathbb{R}^{N \times M}$ ($N > M$). For instance, these nuisance parameters could be some common offsets such as the transmit time, the clock bias, and the transmit power in time-of-arrival (TOA) or received signal strength (RSS) based localization [2], or they could represent some redundant signals like the undesired signatures in hyperspectral imaging [3]. In fact, an estimation problem with linear nuisance parameters widely exists in many other fields such as communications [4–6], source separation [7], and machine learning [8, 9]. Though only

*Correspondence: Y.hu-1@tudelft.nl
The Faculty of Electrical Engineering, Mathematics and Computer Science, Delft University of Technology, Mekelweg 4, 2628 CD Delft, The Netherlands

Bayesian approaches are generally studied in case of nuisance parameters [1, 10, 11], in this paper, we mainly adopt deterministic approaches, for which we first formulate our general model with linear nuisance parameters as

$$\mathbf{y} = \mathbf{H}\mathbf{x} + \mathbf{G}\mathbf{u} + \mathbf{n}, \quad (1)$$

where we assume that

1. The concatenation of $\mathbf{H}$ and $\mathbf{G}$ has full column rank, i.e., Rank([ $\mathbf{H}\ \mathbf{G}$ ]) $= L + M$;
2. The noise $\mathbf{n}$ is zero-mean, i.e., the expected value of $\mathbf{n}$ is $E(\mathbf{n}) = \mathbf{0}$;
3. The noise $\mathbf{n}$ is white (e.g., after whitening), i.e., the covariance matrix $\mathbf{\Sigma_n}$ is (scaled) identity $\mathbf{\Sigma_n} = \sigma^2 \mathbf{I}_N$, where $\mathbf{I}_N$ is the $N \times N$ identity matrix.

Note the noise $\mathbf{n}$ does not have to be Gaussian distributed[1], although it is true for many cases.

To cope with this kind of problem in case $\mathbf{u}$ is deterministic, three methods are often considered: (1) the joint estimation approach estimates the unknown $\mathbf{x}$ together with the unknown nuisance term $\mathbf{u}$ (e.g., the location and the unknown clock bias in [12]); (2) the orthogonal subspace projection (OSP) approach projects out the nuisance term **u** such that the resulting observation vector is only subject to **x** (e.g., the extraction of the desired signature in [13]); (3) the differential signal processing approach firstly chooses a reference and then estimates **x** from the differences between the reference and the observations [14–18]. Note that these methods obviously result in three distinct observation sets with different signal-to-noise ratios (SNRs), which will greatly influence the estimation performance. Therefore, a vast amount of research has been conducted on these methods, though all follow separate paths. Admittedly, some early results have been reported bridging the first two methods. For instance, the famous OSP-based solution using a matched filter to maximize the output SNR proposed in [19] was later on proven to be equivalent to the least squares (LS) approach based on the joint estimation [20, 21]. However, the proposed differential approaches are still widely regarded as a common but distinct way to cope with linear nuisance parameters. One of the most famous applications is time-based localization (TOA or time-difference-of-arrival (TDOA)), where many papers exist on selecting an optimal reference [22–24], constructing an optimal observation subset [25–27] or just using the full observation set adopting each sample as a reference [28–30]. All these issues never occur in the first two methods due to the fact that they are free of a reference. In a nutshell, there still seems to be a huge and inevitable gap between the differential approaches and the other two.

This paper analytically investigates the relations between all three methods, where the corresponding best linear unbiased estimators (BLUEs) are presented and discussed. Since the general framework in (1) is used throughout this paper, all the conclusions apply to any kind of problem that can be written in this form, which is exactly the strength of this paper. We also present some localization examples to verify our conclusions. To summarize, the main contributions of this paper are listed below.

1. For the first time, we extend the differential signal processing approach to a more general framework, which can cope with multiple nuisance parameters, whereas most existing methods consider a single nuisance parameter.
2. Surprisingly, the BLUEs of the three considered methods are proven rigorously to be identical to each other if an appropriate preprocessing step is used. This might be expected or known w.r.t. the first two methods, but the equivalence with differential methods has never been reported before.
3. Compared with the joint estimation method, which directly utilizes all the original observations, none of the other two methods suffers any information loss.
4. Although differential methods seem to rely on the selected reference, selecting the right reference is not important since there is no actual trace of the selected reference in the corresponding BLUE. This is in sharp contrast to what is commonly considered in literature.
5. As far as the differencing process is concerned, the differential observation set associated with a single reference already preserves the full data information.

The rest of this paper is organized as follows. Section 2 presents the relations between the three considered methods. Some examples of source localization are shown and numerically studied to support our conclusions in Section 3. Finally, Section 4 summarizes this paper.

## 2 Handling linear nuisance parameters

In this section, we study the relations between the joint estimation, the OSP-based estimation, and the differential estimation by investigating their corresponding BLUEs, where for the first time, a general differential approach is introduced coping with multiple nuisance parameters.

### 2.1 Joint estimation

The joint least squares (JLS) estimate of **x** and **u**, based on the model (1), is given by

$$\begin{bmatrix} \hat{\mathbf{x}}_{jls} \\ \hat{\mathbf{u}}_{jls} \end{bmatrix} = \left( \begin{bmatrix} \mathbf{H}^T \\ \mathbf{G}^T \end{bmatrix} \begin{bmatrix} \mathbf{H} & \mathbf{G} \end{bmatrix} \right)^{-1} \begin{bmatrix} \mathbf{H}^T \\ \mathbf{G}^T \end{bmatrix} \mathbf{y}, \qquad (2)$$

where we have used the fact that the augmented matrix $\begin{bmatrix} \mathbf{H} & \mathbf{G} \end{bmatrix}$ has a full column rank. Obviously, $\hat{\mathbf{x}}_{jls}$ is the BLUE, since **n** is the zero-mean white noise, according to the Gauss-Markov theorem [1].

### 2.2 OSP-based estimation

If we prefer to project out the nuisance term **u**, an orthogonal subspace projector can be formulated [19] as

$$\mathbf{P}_{\mathbf{G}}^{\perp} \triangleq \mathbf{I}_N - \mathbf{G}\mathbf{G}^{\dagger}, \qquad (3)$$

where $[\cdot]^{\dagger}$ indicates the pseudo-inverse which is given by $\mathbf{G}^{\dagger} \triangleq (\mathbf{G}^T\mathbf{G})^{-1}\mathbf{G}^T$, since **G** is assumed to have a full column rank. Applying $\mathbf{P}_{\mathbf{G}}^{\perp}$ to our original model in (1) results in a new model

$$\mathbf{P}_{\mathbf{G}}^{\perp}\mathbf{y} = \mathbf{P}_{\mathbf{G}}^{\perp}\mathbf{H}\mathbf{x} + \mathbf{P}_{\mathbf{G}}^{\perp}\mathbf{n}, \qquad (4)$$

where the impact of the nuisance term **u** is eliminated. Due to the symmetry and the idempotence of an orthogonal subspace projector, i.e., $\mathbf{P}_{\mathbf{G}}^{\perp} = \mathbf{P}_{\mathbf{G}}^{\perp T}$ and $\mathbf{P}_{\mathbf{G}}^{\perp} = \mathbf{P}_{\mathbf{G}}^{\perp 2}$, we obtain the covariance matrix of the model noise in (4) as $\mathbf{\Sigma}_{\mathbf{P}_{\mathbf{G}}^{\perp}\mathbf{n}} = \sigma^2 \mathbf{P}_{\mathbf{G}}^{\perp} \mathbf{P}_{\mathbf{G}}^{\perp T} = \sigma^2 \mathbf{P}_{\mathbf{G}}^{\perp}$. Then, following the OSP-based model (4), the corresponding LS optimization problem can be formulated as

$$\min_{\mathbf{x}} ||\mathbf{P}_{\mathbf{G}}^{\perp}\mathbf{y} - \mathbf{P}_{\mathbf{G}}^{\perp}\mathbf{H}\mathbf{x}||_2^2, \qquad (5)$$

which leads to the following OSP-based LS estimate Type I of **x**

$$\hat{\mathbf{x}}_{\text{osp}-1} = (\mathbf{H}^T \mathbf{P}_\mathbf{G}^{\perp T} \mathbf{P}_\mathbf{G}^\perp \mathbf{H})^{-1} \mathbf{H}^T \mathbf{P}_\mathbf{G}^{\perp T} \mathbf{P}_\mathbf{G}^\perp \mathbf{y}$$
$$= (\mathbf{H}^T \mathbf{P}_\mathbf{G}^\perp \mathbf{H})^{-1} \mathbf{H}^T \mathbf{P}_\mathbf{G}^\perp \mathbf{y}. \qquad (6)$$

However, the model noise $\mathbf{P}_\mathbf{G}^\perp \mathbf{n}$ in (4) is not white, i.e., $\mathbf{\Sigma}_{\mathbf{P}_\mathbf{G}^\perp \mathbf{n}}$ is not a (scaled) identity. Moreover, the orthogonal subspace projector $\mathbf{P}_\mathbf{G}^\perp$ is obviously singular, which implies that the covariance matrix $\mathbf{\Sigma}_{\mathbf{P}_\mathbf{G}^\perp \mathbf{n}}$ is not invertible and hence can not be used to whiten the model (4). Therefore, it is very difficult to decide at this point whether $\hat{\mathbf{x}}_{\text{osp}-1}$ is the BLUE or not.

To cope with that, we need to introduce another type of OSP-based LS estimator for **x**. If this estimator can be shown to be the BLUE and can also be proven equivalent to $\hat{\mathbf{x}}_{\text{osp}-1}$, then we can conclude that both of them are the BLUE.

Assume that $\mathbf{U}_n \in \mathbb{R}^{N \times (N-M)}$ contains orthonormal basis vectors spanning the null space of **G**. Then, the idea of this second OSP-based estimator is to adopt the null space of **G** to remove the impact of **u**. More specifically, pre-multiplying $\mathbf{U}_n^T$ on both sides of our original model leads to

$$\mathbf{U}_n^T \mathbf{y} = \mathbf{U}_n^T \mathbf{H} \mathbf{x} + \mathbf{U}_n^T \mathbf{n}. \qquad (7)$$

Note that (4) can be obtained from (7) by multiplying it on both sides with $\mathbf{U}_n$ since $\mathbf{U}_n \mathbf{U}_n^T = \mathbf{P}_\mathbf{G}^\perp$ [31], and hence these two models are basically equivalent. We can also see that, since $\mathbf{U}_n$ is an isometry, the model noise $\mathbf{U}_n^T \mathbf{n}$ remains white, i.e., the covariance matrix of $\mathbf{U}_n^T \mathbf{n}$ is $\mathbf{\Sigma}_{\mathbf{U}_n^T \mathbf{n}} = \sigma^2 \mathbf{U}_n^T \mathbf{U}_n = \sigma^2 \mathbf{I}_{N-M}$, which means that the LS estimate of this model is the BLUE.

Applying the LS criterion to the model (7) results in the optimization problem

$$\min_\mathbf{x} || \mathbf{U}_n^T \mathbf{y} - \mathbf{U}_n^T \mathbf{H} \mathbf{x} ||_2^2, \qquad (8)$$

from which we can obtain the OSP-based LS estimate type II of **x** as

$$\hat{\mathbf{x}}_{\text{osp}-2} = (\mathbf{H}^T \mathbf{U}_n \mathbf{U}_n^T \mathbf{H})^{-1} \mathbf{H}^T \mathbf{U}_n \mathbf{U}_n^T \mathbf{y}. \qquad (9)$$

Due to the fact that $\mathbf{U}_n \mathbf{U}_n^T = \mathbf{P}_\mathbf{G}^\perp$, we obtain the equivalence $\hat{\mathbf{x}}_{\text{osp}-1} \equiv \hat{\mathbf{x}}_{\text{osp}-2}$ and hence both estimators represent the BLUE. In the later simulations, these two OSP-based BLUEs will be considered together for convenience.

Finally, to end this subsection, we would like to focus on the equivalence between the joint estimation and the OSP-based estimation approaches. In fact, the equivalence between $\hat{\mathbf{x}}_{\text{jls}}$ and $\hat{\mathbf{x}}_{\text{osp}-1}$ is already known [20, 21, 32], but we found it useful to revisit this result from a different viewpoint. To be explicit, applying the block-wise inversion to (2), we can easily rewrite the joint LS estimate of **x** and **u** as

$$\begin{bmatrix} \hat{\mathbf{x}}_{\text{jls}} \\ \hat{\mathbf{u}}_{\text{jls}} \end{bmatrix} = \begin{bmatrix} \mathbf{M_G} & -\mathbf{M_G} \mathbf{H}^T (\mathbf{G}^\dagger)^T \\ -\mathbf{M_H} \mathbf{G}^T (\mathbf{H}^\dagger)^T & \mathbf{M_H} \end{bmatrix} \begin{bmatrix} \mathbf{H}^T \\ \mathbf{G}^T \end{bmatrix} \mathbf{y},$$
$$= \begin{bmatrix} \mathbf{M_G} \mathbf{H}^T - \mathbf{M_G} \mathbf{H}^T (\mathbf{G}^\dagger)^T \mathbf{G}^T \\ \mathbf{M_H} \mathbf{G}^T - \mathbf{M_H} \mathbf{G}^T (\mathbf{H}^\dagger)^T \mathbf{H}^T \end{bmatrix} \mathbf{y},$$
$$= \begin{bmatrix} \mathbf{M_G} \mathbf{H}^T \mathbf{P}_\mathbf{G}^\perp \\ \mathbf{M_H} \mathbf{G}^T \mathbf{P}_\mathbf{H}^\perp \end{bmatrix} \mathbf{y}, \qquad (10)$$

where $\mathbf{M_G} \triangleq (\mathbf{H}^T \mathbf{P}_\mathbf{G}^\perp \mathbf{H})^{-1}$ and $\mathbf{M_H} \triangleq (\mathbf{G}^T \mathbf{P}_\mathbf{H}^\perp \mathbf{G})^{-1}$ with $\mathbf{P}_\mathbf{H}^\perp \triangleq \mathbf{I} - \mathbf{H} \mathbf{H}^\dagger$. From (10), we can directly observe that $\hat{\mathbf{x}}_{\text{jls}} = \mathbf{M_G} \mathbf{H}^T \mathbf{P}_\mathbf{G}^\perp \mathbf{y}$ and hence

$$\hat{\mathbf{x}}_{\text{jls}} \equiv \hat{\mathbf{x}}_{\text{osp}-1} \equiv \hat{\mathbf{x}}_{\text{osp}-2},$$

where the equivalence between $\hat{\mathbf{x}}_{\text{jls}}$ and $\hat{\mathbf{x}}_{\text{osp}-2}$ is an interesting observation that has never been directly reported before, to the best of our knowledge.

### 2.3 Differential signal processing

In this subsection, we would like to examine differential approaches. This method firstly selects a reference and then removes the impact of **u** by taking differences between the observations and the reference. To be specific, if the $j$th observation $y_j$ is selected as the reference, a new differential observation set can be constructed as

$$\mathbf{d}_j \triangleq \begin{bmatrix} \vdots \\ y_i - y_j \\ \vdots \end{bmatrix}_{(N-1) \times 1} = \mathbf{\Gamma}_j \mathbf{y}, \; i \neq j, \qquad (11)$$

where

$$\mathbf{\Gamma}_j \triangleq \begin{bmatrix} \mathbf{I}_{j-1} & -\mathbf{1}_{(j-1) \times 1} & \mathbf{0} \\ \mathbf{0} & -\mathbf{1}_{(N-j) \times 1} & \mathbf{I}_{N-j} \end{bmatrix}_{(N-1) \times N} \qquad (12)$$

with **1**, the all-one matrix (sizes are mentioned in subscript if needed) and the size of the observation set are reduced to $N - 1$ since $j$ is fixed for every element in $\mathbf{d}_j$. This type of observation set is very popular and has wide applications in source localization and many other areas. Clearly, it can only be used to remove a single nuisance parameter in case $\mathbf{G} = \mathbf{1}_{N \times 1}$.

One may also suggest to select the average of the observations as the reference [16, eq. (28)], thus leading to another kind of differential observation set, given by

$$\mathbf{d}_{\text{avg}} \triangleq \begin{bmatrix} \vdots \\ y_i - \bar{y} \\ \vdots \end{bmatrix}_{N \times 1} = \mathbf{P}_{\mathbf{1}_{N \times 1}}^\perp \mathbf{y} \qquad (13)$$

where $\mathbf{P}_{\mathbf{1}_{N \times 1}}^\perp \triangleq \mathbf{I} - \mathbf{1}_{N \times 1} \mathbf{1}_{N \times 1}^\dagger = \mathbf{I}_N - \frac{1}{N} \mathbf{1}_{N \times N}$. Sometimes, the use of this type of observation set to eliminate the nuisance parameters can be implicit [4], i.e., taking the average of the observations is not clearly pointed out. However, this case can obviously be linked to the OSP-based estimation with a single nuisance parameter in case $\mathbf{G} = \mathbf{1}_{N \times 1}$. Therefore, we are more interested in the simple differencing process of (11), where the reference index $j$ seems to play a significant role.

As already pointed out, (11) only eliminates one nuisance parameter. Nevertheless, we would like to extend this to tackle multiple nuisance parameters, i.e., we would like to relax the constraint $\mathbf{G} = \mathbf{1}_{N \times 1}$ to $\text{rank}(\mathbf{G}) = M \geq 1$. The idea we will adopt here is based on eliminating the impact of the nuisance parameters one by one, which requires $M$ differencing steps.

To achieve that, we write $\mathbf{G} = [\mathbf{g}_1, \cdots, \mathbf{g}_M]$ with $\mathbf{g}_k$ the $k$th column vector of $\mathbf{G}$ related to the $k$th nuisance parameter $u_k$ ($1 \leq k \leq M$). Thus, our original model in (1) can be rewritten as

$$\mathbf{y} = \mathbf{H}\mathbf{x} + \underbrace{\mathbf{g}_1 u_1 + \cdots + \mathbf{g}_M u_M}_{M \text{ nuisance parameters}} + \mathbf{n}. \qquad (14)$$

We then eliminate the nuisance parameters recursively in the order of $u_1, \cdots, u_M$, although the explicit ordering is not important. At the $k$th iteration, when $k-1$ nuisance parameters have already been canceled, the observation vector containing the remaining nuisance parameters can be written as

$$\mathbf{d}^{(k-1)} = \mathbf{H}^{(k-1)}\mathbf{x} + \underbrace{\mathbf{g}_k^{(k-1)} u_k + \cdots + \mathbf{g}_M^{(k-1)} u_M}_{M-k+1 \text{ nuisance parameters}} + \mathbf{n}^{(k-1)}, \qquad (15)$$

where the superscript $(\cdot)^{(k-1)}$ indicates the variables after $k-1$ differencing steps, $\mathbf{y}^{(k-1)}, \mathbf{g}_k^{(k-1)}, \cdots, \mathbf{g}_M^{(k-1)}, \mathbf{n}^{(k-1)} \in \mathbb{R}^{(N-k+1) \times 1}$ and $\mathbf{H}^{(k-1)} \in \mathbb{R}^{(N-k+1) \times L}$. We also assume that, for $k = 1$, $\mathbf{d}^{(0)} = \mathbf{y}$ and similarly $\mathbf{H}^{(0)} = \mathbf{H}, \mathbf{g}_k^{(0)} = \mathbf{g}_k$, and $\mathbf{n}^{(0)} = \mathbf{n}$.

To cancel $u_k$, we first notice that some elements of $\mathbf{g}_k^{(k-1)}$ might be zero, i.e, $u_k$ yields no impact on the corresponding observations in $\mathbf{d}^{(k-1)}$ and hence these observations should not be involved in the differencing process at this iteration. Without loss of generality, we assume that the first $K$ elements of $\mathbf{g}_k^{(k-1)}$ are zero, where $1 \leq K \leq N-k-1$ (there should be at least two non-zero elements for executing the differencing process). Then, among the remaining observations impacted by $u_k$, we select the $j$th element as the reference, $K+1 \leq j \leq N-k+1$, and perform the following differencing step

$$\mathbf{d}^{(k)} = \begin{bmatrix} [\mathbf{d}^{(k-1)}]_1 \\ \vdots \\ [\mathbf{d}^{(k-1)}]_K \\ \frac{[\mathbf{d}^{(k-1)}]_{K+1}}{[\mathbf{g}_k^{(k-1)}]_{K+1}} - \frac{[\mathbf{d}^{(k-1)}]_j}{[\mathbf{g}_k^{(k-1)}]_j} \\ \vdots \\ \frac{[\mathbf{d}^{(k-1)}]_{j-1}}{[\mathbf{g}_k^{(k-1)}]_{j-1}} - \frac{[\mathbf{d}^{(k-1)}]_j}{[\mathbf{g}_k^{(k-1)}]_j} \\ \frac{[\mathbf{d}^{(k-1)}]_{j+1}}{[\mathbf{g}_k^{(k-1)}]_{j+1}} - \frac{[\mathbf{d}^{(k-1)}]_j}{[\mathbf{g}_k^{(k-1)}]_j} \\ \vdots \\ \frac{[\mathbf{d}^{(k-1)}]_{N-k+1}}{[\mathbf{g}_k^{(k-1)}]_{N-k+1}} - \frac{[\mathbf{d}^{(k-1)}]_j}{[\mathbf{g}_k^{(k-1)}]_j} \end{bmatrix}_{(N-k) \times 1} = \mathbf{\Gamma}^{(k)} \mathbf{d}^{(k-1)}, \qquad (16)$$

where $\mathbf{\Gamma}^{(k)} \triangleq \begin{bmatrix} \mathbf{I}_K & \mathbf{0} \\ \mathbf{0} & \mathbf{\Gamma}_\perp^{(k)} \text{diag}\left( \left[ \frac{1}{[\mathbf{g}_k^{(k-1)}]_{K+1}}, \cdots, \frac{1}{[\mathbf{g}_k^{(k-1)}]_{N-k-1}} \right]^T \right) \end{bmatrix}$ is the $(N-k) \times (N-k+1)$ differencing operator for $\mathbf{d}^{(k-1)}$ with

$$\mathbf{\Gamma}_\perp^{(k)} \triangleq \begin{bmatrix} \mathbf{I}_{j-K-1} & -\mathbf{1}_{(j-K-1) \times 1} & \mathbf{0} \\ \mathbf{0} & -\mathbf{1}_{(N-k-j+1) \times 1} & \mathbf{I}_{N-k-j+1} \end{bmatrix}, \qquad (17)$$

and obviously $\mathbf{\Gamma}^{(k)} \mathbf{g}_k^{(k-1)} = \mathbf{0}$. Accordingly, the new differential observation vector $\mathbf{d}^{(k)}$ can be formulated as

$$\mathbf{d}^{(k)} = \underbrace{\mathbf{\Gamma}^{(k)} \mathbf{H}^{(k-1)}}_{\mathbf{H}^{(k)}} \mathbf{x} + \underbrace{\mathbf{\Gamma}^{(k)} \mathbf{g}_{k+1}^{(k-1)}}_{\mathbf{g}_{k+1}^{(k)}} u_{k+1} + \cdots + \underbrace{\mathbf{\Gamma}^{(k)} \mathbf{g}_M^{(k-1)}}_{\mathbf{g}_M^{(k)}} u_M + \underbrace{\mathbf{\Gamma}^{(k)} \mathbf{n}^{(k-1)}}_{\mathbf{n}^{(k)}},$$
$$\underbrace{\phantom{XX}}_{M-k \text{ nuisance parameters}} \qquad (18)$$

where $u_k$ has been canceled.

We can see that (18) is similar to (15) with $k-1$ replaced by $k$. So it is clear that this recursive process can remove all nuisance parameters. Note that the number of zero values $K$ as well as the reference index $j$ could be different in every step, but for simplicity, we use the same notation in every step.

To understand the interaction of the successive differencing steps, let us introduce the total differencing operator $\mathbf{\Gamma} = \mathbf{\Gamma}^{(M)} \cdots \mathbf{\Gamma}^{(1)}$, where obviously $\text{rank}(\mathbf{\Gamma}^{(k)} \mathbf{\Gamma}^{(k-1)}) = \text{rank}(\mathbf{\Gamma}^{(k)}) = N-k$ and hence $\mathbf{\Gamma}$ has full row rank. Since it is clear that $\mathbf{\Gamma}\mathbf{G} = \mathbf{0}$, the final differential observation vector $\mathbf{d}^{(M)}$ can be expressed as

$$\mathbf{d}^{(M)} = \mathbf{\Gamma}\mathbf{y} = \mathbf{\Gamma}\mathbf{H}\mathbf{x} + \mathbf{\Gamma}\mathbf{n}, \qquad (19)$$

where the covariance matrix of $\mathbf{\Gamma}\mathbf{n}$ is $\mathbf{\Sigma}_{\mathbf{\Gamma}\mathbf{n}} = \sigma^2 \mathbf{\Gamma}\mathbf{\Gamma}^T$.

Observe that the model noise has become correlated ever since the first step of the differencing process. Therefore, we need to whiten the model in (19) as

$$\begin{aligned}
\boldsymbol{\Sigma}_{\boldsymbol{\Gamma n}}^{-1/2}\mathbf{d}^{(M)} &= \boldsymbol{\Sigma}_{\boldsymbol{\Gamma n}}^{-1/2}\boldsymbol{\Gamma}\mathbf{H}\mathbf{x} + \boldsymbol{\Sigma}_{\boldsymbol{\Gamma n}}^{-1/2}\boldsymbol{\Gamma}\mathbf{n}, \\
\Longrightarrow (\boldsymbol{\Gamma}\boldsymbol{\Gamma}^T)^{-1/2}\mathbf{d}^{(M)} &= (\boldsymbol{\Gamma}\boldsymbol{\Gamma}^T)^{-1/2}\boldsymbol{\Gamma}\mathbf{H}\mathbf{x} + (\boldsymbol{\Gamma}\boldsymbol{\Gamma}^T)^{-1/2}\boldsymbol{\Gamma}\mathbf{n}, \\
\Longrightarrow \mathbf{P}\mathbf{y} &= \mathbf{P}\mathbf{H}\mathbf{x} + \mathbf{P}\mathbf{n},
\end{aligned} \quad (20)$$

where the unknown $\sigma^2$ is canceled out at both sides of the equation and $\mathbf{P} \triangleq (\boldsymbol{\Gamma}\boldsymbol{\Gamma}^T)^{-1/2}\boldsymbol{\Gamma}$ which exists since $\boldsymbol{\Gamma}$ has full row rank. Note that $\mathbf{P}$, as well as $\boldsymbol{\Gamma}$ and $\mathbf{d}^{(k)}$, depend on the reference indices $j$ that have been chosen in the successive differencing steps, although this has not been explicitly stated.

Applying the LS criterion, the corresponding optimization problem is now obtained as

$$\min_{\mathbf{x}} ||\mathbf{P}\mathbf{y} - \mathbf{P}\mathbf{H}\mathbf{x}||_2^2, \quad (21)$$

which leads to the following BLUE for model (19)

$$\hat{\mathbf{x}}_d = (\mathbf{H}^T\mathbf{P}^T\mathbf{P}\mathbf{H})^{-1}\mathbf{H}^T\mathbf{P}^T\mathbf{P}\mathbf{y}. \quad (22)$$

Finally, to prove the equivalence of the estimate $\hat{\mathbf{x}}_d$ to the previous estimates, i.e., to prove that

$$\hat{\mathbf{x}}_{\text{jls}} \equiv \hat{\mathbf{x}}_{\text{osp}-1} \equiv \hat{\mathbf{x}}_{\text{osp}-2} \equiv \hat{\mathbf{x}}_d,$$

we need to establish the relation $\mathbf{P}^T\mathbf{P} = \mathbf{U}_n\mathbf{U}_n^T = \mathbf{P}_{\mathbf{G}}^\perp$. To do that, we first recall that $\boldsymbol{\Gamma}\mathbf{G} = \mathbf{0}$ and that $\boldsymbol{\Gamma}$ has full row rank. Hence, $\boldsymbol{\Gamma}$ can always be written as $\boldsymbol{\Gamma} = \mathbf{Q}\mathbf{U}_n^T$, where $\mathbf{Q}$ is an $(N-M) \times (N-M)$ invertible matrix and $\mathbf{U}_n$ has already been defined before as a basis that spans the null space of $\mathbf{G}$. The proof is completed by computing

$$\begin{aligned}
\mathbf{P}^T\mathbf{P} &= \boldsymbol{\Gamma}^T(\boldsymbol{\Gamma}\boldsymbol{\Gamma}^T)^{-1}\boldsymbol{\Gamma} \\
&= \mathbf{U}_n\mathbf{Q}^T(\mathbf{Q}\mathbf{U}_n^T\mathbf{U}_n\mathbf{Q}^T)^{-1}\mathbf{Q}\mathbf{U}_n^T \\
&= \mathbf{U}_n\mathbf{Q}^T(\mathbf{Q}^T)^{-1}(\mathbf{U}_n^T\mathbf{U}_n)^{-1}\mathbf{Q}^{-1}\mathbf{Q}\mathbf{U}_n^T \\
&= \mathbf{U}_n\mathbf{U}_n^T = \mathbf{P}_{\mathbf{G}}^\perp,
\end{aligned} \quad (23)$$

where we surprisingly notice that, even though $\mathbf{P}$ and $\boldsymbol{\Gamma}$ are subject to possibly different reference indices $j$, there is no trace of any $j$ in $\mathbf{P}^T\mathbf{P}$ and hence in $\hat{\mathbf{x}}_d$.

**A Simple Illustrative Case:** We would like to demonstrate these three different methods, particularly the differential signal processing, with a simple example. Given $N = 3$ samples, we only assume a single parameter of interest ($L = 1$), but with two linear nuisance parameters ($M = 2$). We also know that $\mathbf{H} = \begin{bmatrix} 3 & 6 & 7 \end{bmatrix}^T$ and $\mathbf{G} = \begin{bmatrix} 3 & 5 & 2 \\ 2 & 4 & 8 \end{bmatrix}^T$ and hence the joint estimator in (2) results into $\begin{bmatrix} \hat{x}_{\text{jls}} \\ \hat{\mathbf{u}}_{\text{jls}} \end{bmatrix} = \begin{bmatrix} -3.2 & 2 & -0.2 \\ 2 & -1 & 0 \\ 2.3 & -1.5 & 0.3 \end{bmatrix}\mathbf{y}$, where the parameter estimate of interest is given by $\hat{x}_{\text{jls}} = \begin{bmatrix} -3.2 & 2 & -0.2 \end{bmatrix}\mathbf{y}$. Then, we calculate $\mathbf{P}_{\mathbf{G}}^\perp = \begin{bmatrix} 0.7171 & -0.4482 & 0.0448 \\ -0.4482 & 0.2801 & -0.0280 \\ 0.0448 & -0.0280 & 0.0028 \end{bmatrix}$ and $\mathbf{U}_n = \begin{bmatrix} -0.8468 & 0.5293 & -0.0529 \end{bmatrix}^T$ such that two OSP-based estimators in (6) and (9) can easily be carried out and proved to be equal to $\hat{x}_{\text{jls}}$. We will not present more details for simplicity but particularly focus on the differential method. Since there exist two linear nuisance parameters, it would take two steps for eliminating all of them:

1. In the first step ($k = 1$), we arbitrarily select the third element of $\mathbf{y}$ as the reference ($j = 3$). Splitting $\mathbf{G}$ by columns, we have $\mathbf{g}_1^{(0)} = [3\ 5\ 2]^T\ \mathbf{g}_2^{(0)} = [2\ 4\ 8]^T$. According to (16), the new differential observation vector can be obtained as $\mathbf{d}^{(1)} = \begin{bmatrix} y_1/3 - y_3/2 & y_2/5 - y_3/2 \end{bmatrix}^T = \boldsymbol{\Gamma}^{(1)}\mathbf{y}$, where $\boldsymbol{\Gamma}^{(1)} = \begin{bmatrix} 1/3 & 0 & -1/2 \\ 0 & 1/5 & -1/2 \end{bmatrix}$. We can observe from $\boldsymbol{\Gamma}^{(1)}\mathbf{G} = \begin{bmatrix} 0 & -10/3 \\ 0 & -16/5 \end{bmatrix}$ that the impact of the first nuisance parameter $u_1$ is already eliminated. Also, $\mathbf{g}_2^{(1)} = \boldsymbol{\Gamma}^{(1)}\mathbf{g}_2^{(0)}$ corresponds to the last column and the next nuisance parameter $u_2$.
2. In the second step ($k = 2$), the first element of $\mathbf{d}^{(1)}$ is selected as the reference ($j = 1$). The differential observation becomes a scalar as $\mathbf{d}^{(2)} = -\frac{5}{16}(y_2/5 - y_3/2) + \frac{3}{10}(y_1/3 - y_3/2) = \boldsymbol{\Gamma}^{(2)}\mathbf{d}^{(1)} = \boldsymbol{\Gamma}\mathbf{y}$, where $\boldsymbol{\Gamma}^{(2)} = [3/10\ -5/16]$ and $\boldsymbol{\Gamma} = \boldsymbol{\Gamma}^{(2)}\boldsymbol{\Gamma}^{(1)} = [1/10\ -1/16\ 1/160]$. Now, we can readily observe that all the nuisance parameters are eliminated, since $\boldsymbol{\Gamma}\mathbf{G} = \mathbf{0}$.

With a known $\boldsymbol{\Gamma}$, we can easily whiten the model in (19) and obtain the differential estimator in (22). Moreover, the equivalence of the differential estimation can also be proved by observing $\mathbf{P}^T\mathbf{P} = \boldsymbol{\Gamma}^T(\boldsymbol{\Gamma}\boldsymbol{\Gamma}^T)^{-1}\boldsymbol{\Gamma} = \mathbf{P}_{\mathbf{G}}^\perp$.

### 2.4 Discussion

We have studied estimation problems in the presence of deterministic linear nuisance parameters based on a general model. Therefore, all the conclusions drawn in this paper are applicable to any optimization problem with a data model that matches our general model (1). The equivalences between the BLUEs of the joint estimation, the OSP-based estimation and the differential estimation are summarized in Table 1 and also in Fig. 1. Some interesting observations from these equivalences are listed below:

1. The joint estimation has to estimate both $\mathbf{x}$ and the nuisance term $\mathbf{u}$ while the other two estimation approaches remove the impact of $\mathbf{u}$ before estimating $\mathbf{x}$.

**Table 1** Relations between the BLUEs related to the joint estimation, the OSP-based estimation, and the differential estimation

| Models | BLUEs | Equality conditions |
|---|---|---|
| Original Model in (1) | Joint estimator in (2) | $[\mathbf{I}_L\ \mathbf{0}_{L\times M}]$[a], $\mathbf{P}_\mathbf{G}^\perp \triangleq \mathbf{I}_N - \mathbf{G}\mathbf{G}^\dagger$ |
| OSP Model Type I in (4) | OSP estimator type I in (6) | $\mathbf{P}_\mathbf{G}^{\perp T}\mathbf{P}_\mathbf{G}^\perp = \mathbf{P}_\mathbf{G}^\perp$ |
| OSP Model Type II in (7) | OSP estimator type II in (9) | $\mathbf{U}_n\mathbf{U}_n^T = \mathbf{P}_\mathbf{G}^\perp$ |
| Differential Model in (19) or the Whitened One in (20) | Differential estimator in (22) | $\mathbf{P}^T\mathbf{P} = \mathbf{P}_\mathbf{G}^\perp$ |

[a] $[\mathbf{I}_L\ \mathbf{0}_{L\times M}]$ is used for extracting $\hat{\mathbf{x}}_{jls}$ in (22)

2. For the OSP-based estimation, in order to remove the impact of **u**, using $\mathbf{P}_\mathbf{G}^\perp$ actually colors the noise, but using $\mathbf{U}_n^T$ keeps the model noise white. Interestingly though, the corresponding LS estimates for those two models are theoretically equivalent and hence they are both the BLUE.

3. In many applications, the differential processing is commonly considered as a separate and independent approach. But, in this paper, we have generally proven its equivalence to the joint estimation and the OSP-based estimation. The differential approach removes the impact of the nuisance parameters by taking differences between the reference and the observations. If one of the observations is selected as a reference, the obtained differential observation set has to be properly whitened in order to obtain the BLUE for this model.

4. From an information theoretic perspective, the joint estimation, which directly utilizes the observations **y**, preserves the full data information, and any preprocessing on the observations might cause an information loss. However, in this paper, all the other BLUEs have been proven to be equivalent to the BLUE of the joint estimation, which implies that neither the OSP-based estimation nor the differential estimation suffers any information loss by removing the impact of the nuisance parameters.

5. It is also worth noting that, for the differential approach, selecting which observation will function as a reference is not important, since the reference index $j$ yields no impact on the BLUE. This is in sharp contrast to what is commonly considered in literature.

6. One might notice that, in the differencing process, $N$ observations can generate a maximum of $N(N-1)/2$ distinct observation differences. In contrast, we only study the estimation problem based on a subset, which is associated with a single reference and corresponds to $N-1$ observation differences. However, from the above conclusions, it is clear that the considered subset already preserves all the information (independent of the reference), which implies that the full set obtains no more information than any subset does. Also this is a novel observation.

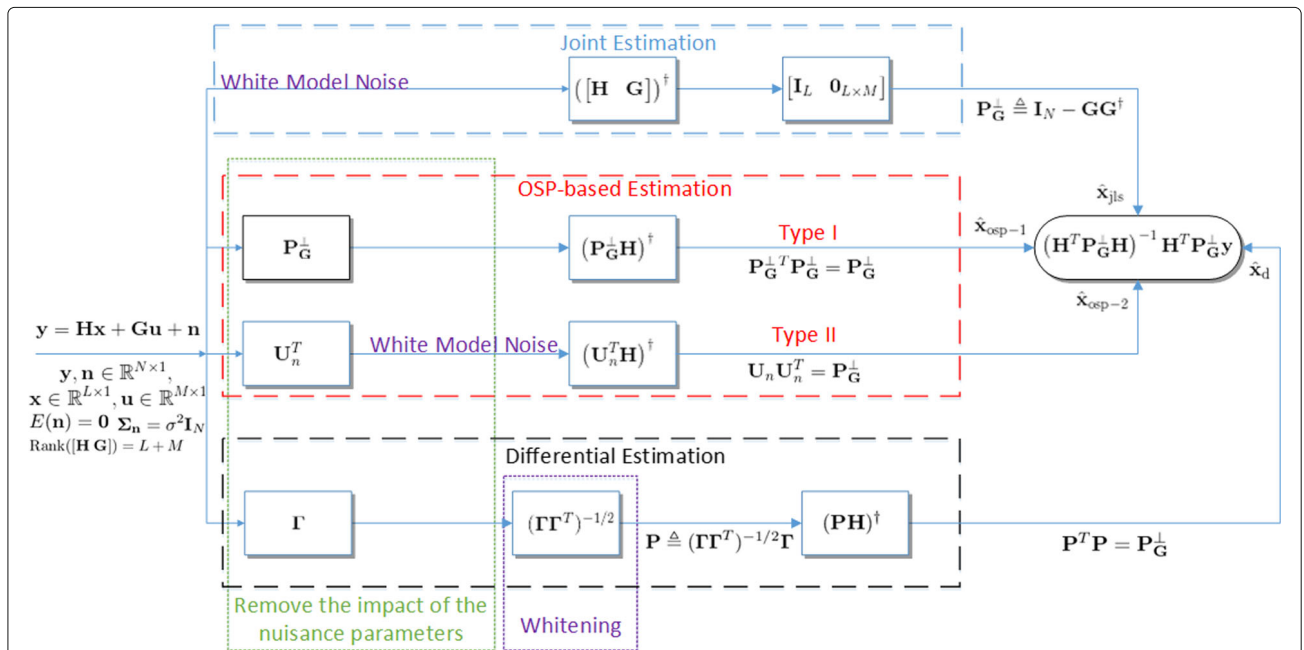

**Fig. 1** Diagram to illustrate the relations between the BLUEs related to the joint estimation, the OSP-based estimation, and the differential estimation. Note that the noise **n** is not necessarily Gaussian distributed and the operator $[\mathbf{I}_L\ \mathbf{0}_{L\times M}]$ is used to extract the first $L$ elements of a vector, i.e., $\hat{\mathbf{x}}_{jls}$

## 3 Localization examples
By studying the relations between the BLUE of the joint estimation, the OSP-based estimation, and the differential estimation, the essence of this paper is to provide some in-depth understanding of coping with unknown nuisance parameters. Some important underlying equivalences have been unveiled, especially the one related to the differential method, since, in many applications, this approach is still considered as a separate optimization problem. Owing to the generality of this paper, one may easily apply our analyses and conclusions to some particular applications, if the data model can be (re)formulated to match our general model (1). Some specific localization examples are detailed next.

### 3.1 Time-based localization
Both TOA- and TDOA-based localization are called time-based localization [2], since they both rely on time measurements (either the global time or the local time). The essence of this kind of localization problem is how to accurately extract distance-related information (e.g., the time of flight (TOF)). Directly using TOA measurements requires not only perfect clock synchronization between the emitters and the receivers but also the knowledge of the transmitting time [33]. In cooperative networks, where clock synchronization is frequently carried out (because the inner clock might drift over time) and the transmitting times are also piggybacked with the transmitted signals, one can precisely calculate the TOFs from the TOA measurements and then localize the target node. However, it is often very expensive to meet those requirements, and most networks are constrained by limited resources and capabilities. Therefore, in most cases, sensors suffer from two linear nuisance parameters, i.e., the unknown clock biases to the global time and the unknown transmitting times.

In this example, we assume $N$ anchor nodes that are perfectly synchronized with the global time and there exists only a clock bias in the target node, which broadcasts beacon signals at unknown local transmit times. We denote $\mathbf{x}_t \in \mathbb{R}^d$ as the target location and $\mathbf{s}_i \in \mathbb{R}^d$ as the $i$th anchor location. For convenience, a single unknown global transmit time $t_0$ is considered for the target node, instead of the local transmit time plus the clock bias. Taking the speed of light $c$ into account, we obtain the TOA measurements as

$$\mathbf{d} = \mathbf{r}(\mathbf{x}_t) + \mathbf{1}_{N \times 1} r_o + \mathbf{n}, \quad (24)$$

where the element $d_i$ of $\mathbf{d}$ indicates the TOA measurement from the $i$th anchor, $\mathbf{r}(\mathbf{x}_t)$ stacks $r_i \triangleq ||\mathbf{x}_t - \mathbf{s}_i||_2$, $r_o \triangleq ct_0$ and $\mathbf{n}$ is the vector of the measurement noise $n_i$ with $\mathbf{n} \sim \mathcal{N}(\mathbf{0}, \sigma^2 \mathbf{I}_N)$. Note that, compared with more realistic scenarios, the model (24) is simplified for convenience but still adequate to make our point.

#### 3.1.1 Taylor series expansion
Obviously, the non-linearity of (24) is a very serious issue for localization, other than the nuisance parameter. Many methods, especially those considering mobile scenarios, directly linearize (24) by a Taylor series expansion (TSE) [34]. Note that this kind of method is very similar to the Gauss-Newton (GN) method [35] and holds the maximum likelihood (ML) property. Since we can obtain the estimate of $\mathbf{x}_t$ by iteratively updating the previous iteration, we first have to apply the TSE to (24) around the target location estimate $\hat{\mathbf{x}}_t^{(k-1)}$ at the $(k-1)$th iteration, thus resulting into

$$\mathbf{d} = \mathbf{r}(\hat{\mathbf{x}}_t^{(k-1)}) + \left.\frac{\partial \mathbf{r}}{\partial \mathbf{x}_t^T}\right|_{\mathbf{x}_t = \hat{\mathbf{x}}_t^{(k-1)}} (\mathbf{x}_t - \hat{\mathbf{x}}_t^{(k-1)}) + \mathbf{1}_{N \times 1} r_o + \mathbf{n}.$$

Then, we rearrange the above equation and present the TSE model for iteration step $k$ as

$$\begin{aligned}
&\mathbf{d} - \mathbf{r}(\hat{\mathbf{x}}_t^{(k-1)}) + \left.\frac{\partial \mathbf{r}}{\partial \mathbf{x}_t^T}\right|_{\mathbf{x}_t = \hat{\mathbf{x}}_t^{(k-1)}} \hat{\mathbf{x}}_t^{(k-1)} \\
&= \left.\frac{\partial \mathbf{r}}{\partial \mathbf{x}_t^T}\right|_{\mathbf{x}_t = \hat{\mathbf{x}}_t^{(k-1)}} \mathbf{x}_t + \mathbf{1}_{N \times 1} r_o + \mathbf{n} \quad (25)\\
&\Rightarrow \boldsymbol{\delta}^{(k-1)} = \boldsymbol{\Delta}^{(k-1)} \mathbf{x}_t + \mathbf{1}_{N \times 1} r_o + \mathbf{n},
\end{aligned}$$

where $\boldsymbol{\delta}^{(k-1)} \triangleq \mathbf{d} - \mathbf{r}(\hat{\mathbf{x}}_t^{(k-1)}) + \boldsymbol{\Delta}^{(k-1)} \hat{\mathbf{x}}_t^{(k-1)}$ and

$$\boldsymbol{\Delta}^{(k-1)} \triangleq \left.\frac{\partial \mathbf{r}}{\partial \mathbf{x}_t^T}\right|_{\mathbf{x}_t = \hat{\mathbf{x}}_t^{(k-1)}} = \left[\ldots, \frac{(\hat{\mathbf{x}}_t^{(k-1)} - \mathbf{s}_i)^T}{||\hat{\mathbf{x}}_t^{(k-1)} - \mathbf{s}_i||_2}, \ldots\right]^T.$$

The localization problem at the $k$th iteration boils down to estimating $\mathbf{x}_t$ from (25) to update the location estimate from the $(k-1)$th iteration. The relation between the TSE model and the general model (1) is presented in Table 2. Note that since the discussed approaches can directly be applied to the TOA measurements with a single nuisance parameter ($M = 1$), the differential approach applied to the TOA measurements actually corresponds to working with the TDOA measurements, i.e.,

$$\begin{bmatrix} \vdots \\ d_{i,j} \\ \vdots \end{bmatrix} = \begin{bmatrix} \vdots \\ d_i - d_j \\ \vdots \end{bmatrix}_{(N-1) \times 1}, i \neq j. \quad (26)$$

However, to avoid any confusion with the TDOA methods we will discuss later on, we will refer to this method as the differential approach applied to the TSE model of the TOA measurements.

**Table 2** Relations between the general model (1) and the considered time-based and RSS-based localization models[a]

| General model (1) | y | H | x | G | u |
|---|---|---|---|---|---|
| TSE model (25) | $\delta^{(k-1)}$ | $\Delta^{(k-1)}$ | $\mathbf{x}_t$ | $\mathbf{1}_{N \times 1}$ | $r_0$ |
| SD-TOA model (31b) | $\mathbf{D}_1' \mathbf{z}_1$ | $\mathbf{D}_1' \mathbf{A}_1'$ | $\mathbf{x}_t$ | $\mathbf{D}_1' \mathbf{A}_1''$ | $[\|\mathbf{x}_t\|_2^2 - r_0^2, r_0]^T$ |
| SD-TDOA model (36b) | $\mathbf{D}_2' \mathbf{z}_2$ | $\mathbf{D}_2' \mathbf{A}_2'$ | $\mathbf{x}_t$ | $\mathbf{D}_2' \mathbf{A}_2''$ | $r_j$ |
| SD-RSS model (43b) | $\mathbf{Dh}$ | $\mathbf{F}'$ | $[\mathbf{x}_t^T, \|\mathbf{x}_t\|_2^2]^T$ | $\mathbf{1}_{N \times 1}$ | $P_0'$ |

[a]All the considered models must be white or whitened, i.e., the covariance of the model noise should be a (scaled) identity

### 3.1.2 Squared distance

The TSE method highly relies on an appropriate initialization that is near the global solution; otherwise, it might converge to a local minimum. Thus, some closed-form solutions were proposed to solve this non-convex problem, which requires squaring the distance norm (SD) for linearization [36]. Unlike the TSE method, the SD method depends on the type of measurements, since different modeling steps are carried out for TOA and TDOA measurements.

**TOA:** Let us first focus on the SD method based on the TOA measurements which can be expressed as

$$d_i = \|\mathbf{x}_t - \mathbf{s}_i\|_2 + r_0 + n_i. \qquad (27)$$

Moving $r_0$ to the other side and squaring both sides of the equation, we obtain

$$(d_i - r_o)^2 = (\|\mathbf{x}_t - \mathbf{s}_i\|_2 + n_i)^2$$
$$\Rightarrow -2\mathbf{s}_i^T \mathbf{x}_t + \|\mathbf{x}_t\|_2^2 - r_0^2 - 2d_i r_0 = d_i^2 - \|\mathbf{s}_i\|_2^2 - 2r_i n_i - n_i^2, \qquad (28)$$

where $r_0^2$ is viewed as a new nuisance parameter. As a result, a linear model with two nuisance parameters ($M = 2$) can be formulated as

$$\mathbf{z}_1 = \mathbf{A}_1 \boldsymbol{\theta}_1 + \boldsymbol{\epsilon}_1, \qquad (29)$$

where $\mathbf{A}_1 \triangleq \begin{bmatrix} \vdots & \vdots & \vdots \\ -2\mathbf{s}_i^T & 1 & -2d_i \\ \vdots & \vdots & \vdots \end{bmatrix}, \boldsymbol{\theta}_1 \triangleq \begin{bmatrix} \mathbf{x}_t \\ \|\mathbf{x}_t\|_2^2 - r_0^2 \\ r_0 \end{bmatrix}, \mathbf{z}_1 \triangleq \begin{bmatrix} \vdots \\ d_i^2 - \|\mathbf{s}_i\|_2^2 \\ \vdots \end{bmatrix}$ and

$$\boldsymbol{\epsilon}_1 \triangleq \begin{bmatrix} \vdots \\ 2r_i n_i + n_i^2 \\ \vdots \end{bmatrix} \approx \begin{bmatrix} \vdots \\ 2r_i n_i \\ \vdots \end{bmatrix} = 2\mathbf{D}_1 \mathbf{n}. \qquad (30a)$$

Here, we denote $\mathbf{D}_1 = \text{diag}([r_1, \cdots, r_N]^T)$ with $\text{diag}(\cdot)$ as a diagonal matrix with its argument on the diagonal, and hence $\boldsymbol{\Sigma}_{\boldsymbol{\epsilon}_1} = 4\sigma^2 \mathbf{D}_1^2$. This SD-TOA model is widely considered [37–41]. Some researchers apply the differencing process to remove the nuisance parameters [24, 33, 42–45] while some others use the OSP method [16, 46]. Note that the model noise in (30a) is still not white, and hence, an appropriate whitening procedure is required. Assuming $\mathbf{D}_1$ is perfectly known, we can whiten the model (29) as

$$\boldsymbol{\Sigma}_{\boldsymbol{\epsilon}_1}^{-1/2} \mathbf{z}_1 = \boldsymbol{\Sigma}_{\boldsymbol{\epsilon}_1}^{-1/2} \mathbf{A}_1 \boldsymbol{\theta}_1 + \boldsymbol{\Sigma}_{\boldsymbol{\epsilon}_1}^{-1/2} \boldsymbol{\epsilon}_1 \qquad (31a)$$
$$\Rightarrow \mathbf{D}_1' \mathbf{z}_1 = \mathbf{D}_1' \mathbf{A}_1 \boldsymbol{\theta}_1 + \mathbf{D}_1' \boldsymbol{\epsilon}_1 \qquad (31b)$$

where $\mathbf{D}_1' \triangleq \mathbf{D}_1^{-1}$ and the covariance matrix of $\mathbf{D}_1' \boldsymbol{\epsilon}_1$ is now a scaled identity, i.e., $\boldsymbol{\Sigma}_{\mathbf{D}_1' \boldsymbol{\epsilon}_1} = 4\sigma^2 \mathbf{I}_N$. In practice, a LS estimate based on the model (29) can first be used to construct an estimate of $\mathbf{D}_1$ for carrying out the whitening. Then, the estimate of $\mathbf{D}_1$ can be repeatedly updated to approach the true $\mathbf{D}_1$ with a more accurate location estimate. In this paper though, we only want to evaluate its best performance and hence directly use the true $\mathbf{D}_1$. Finally, expressing $\mathbf{A}_1 = [\mathbf{A}_1', \mathbf{A}_1'']$ with $\mathbf{A}_1'$ and $\mathbf{A}_1''$, respectively, containing the first $d$ and the remaining columns, the relation between the whitened SD-TOA model and the general model (1) is presented in Table 2.

**TDOA:** Directly applying the differencing process on the TOA observations $\mathbf{d}$ removes the unknown nuisance parameter $r_0$, resulting in the TDOA measurements

$$d_{i,j} = \|\mathbf{x}_t - \mathbf{s}_i\|_2 - \|\mathbf{x}_t - \mathbf{s}_j\|_2 + n_{i,j}, i \neq j, \qquad (32)$$

where $n_{i,j} = n_i - n_j$. Introducing $r_j = \|\mathbf{x}_t - \mathbf{s}_j\|_2$ as a new unknown parameter, we can linearize (32) using the following squaring operation

$$(d_{i,j} + r_j)^2 = (\|\mathbf{x}_t - \mathbf{s}_j - (\mathbf{s}_i - \mathbf{s}_j)\|_2 + n_{i,j})^2$$
$$\Rightarrow -2(\mathbf{s}_i - \mathbf{s}_j)^T \mathbf{x}_t - 2d_{i,j} r_j = d_{i,j}^2 + \|\mathbf{s}_j\|_2^2 - \|\mathbf{s}_i\|_2^2 - 2r_i n_{i,j} - n_{i,j}^2. \qquad (33)$$

As a result, a linear model with a single unknown nuisance parameter $r_j$ ($M = 1$) can be formulated as

$$\mathbf{z}_2 = \mathbf{A}_2 \boldsymbol{\theta}_2 + \boldsymbol{\epsilon}_2, \qquad (34)$$

where $\mathbf{A}_2 \triangleq -2 \begin{bmatrix} \vdots & \vdots \\ (\mathbf{s}_i - \mathbf{s}_j)^T & d_{i,j} \\ \vdots & \vdots \end{bmatrix}, \boldsymbol{\theta}_2 \triangleq \begin{bmatrix} \mathbf{x}_t \\ r_j \end{bmatrix}, \mathbf{z}_2 \triangleq$

$\begin{bmatrix} \vdots \\ d_{i,j}^2 + ||\mathbf{s}_j||_2^2 - ||\mathbf{s}_i||_2^2 \\ \vdots \end{bmatrix}$, and

$$\boldsymbol{\epsilon}_2 \triangleq \begin{bmatrix} \vdots \\ 2r_i n_{i,j} + n_{i,j}^2 \\ \vdots \end{bmatrix} \approx \begin{bmatrix} \vdots \\ 2r_i n_{i,j} \\ \vdots \end{bmatrix} = 2\mathbf{D}_2 \boldsymbol{\Gamma}_j \mathbf{n}. \quad (35a)$$

Here, we denote $\mathbf{D}_2 = \mathrm{diag}([\cdots, r_i, \cdots]^T), i \neq j$, and hence, $\boldsymbol{\Sigma}_{\boldsymbol{\epsilon}_2} = 4\sigma^2 \mathbf{D}_2 \boldsymbol{\Gamma}_j \boldsymbol{\Gamma}_j^T \mathbf{D}_2^T$. Also, this SD-TDOA model has been commonly adopted in literature [14, 33, 47–51]. Among the TDOA localization techniques based on this model, the famous Chan algorithm [14], from which many others stem, is actually equivalent to some earlier works [52–54], where the unknown $r_j$ is simply removed by the OSP method. Again, note that the model noise (35a) is not white. Assuming $\mathbf{D}_2$ is perfectly known (as already explained for $\mathbf{D}_1$, in practice, $\mathbf{D}_2$ should be iteratively estimated), we can whiten the model (34) as

$$\boldsymbol{\Sigma}_{\boldsymbol{\epsilon}_2}^{-1/2} \mathbf{z}_2 = \boldsymbol{\Sigma}_{\boldsymbol{\epsilon}_2}^{-1/2} \mathbf{A}_2 \boldsymbol{\theta}_2 + \boldsymbol{\Sigma}_{\boldsymbol{\epsilon}_2}^{-1/2} \boldsymbol{\epsilon}_2 \quad (36a)$$
$$\Rightarrow \mathbf{D}_2' \mathbf{z}_2 = \mathbf{D}_2' \mathbf{A}_2 \boldsymbol{\theta}_2 + \mathbf{D}_2' \boldsymbol{\epsilon}_2, \quad (36b)$$

where $\mathbf{D}_2' \triangleq (\mathbf{D}_2 \boldsymbol{\Gamma}_j \boldsymbol{\Gamma}_j^T \mathbf{D}_2^T)^{-1/2}$ and the covariance matrix of $\mathbf{D}_2' \boldsymbol{\epsilon}_2$ is now a scaled identity, i.e., $\boldsymbol{\Sigma}_{\mathbf{D}_2' \boldsymbol{\epsilon}_2} = 4\sigma^2 \mathbf{I}_{N-1}$. Finally, we split $\mathbf{A}_2$ into $\mathbf{A}_2 = [\mathbf{A}_2', \mathbf{A}_2'']$ with $\mathbf{A}_2'$ and $\mathbf{A}_2''$, respectively, containing the first $d$ and the remaining columns. The relation between the whitened SD-TDOA model and the general model (1) is finally presented in Table 2.

**Numerical results:** We have conducted a Monte Carlo simulation with 1000 trials to verify our conclusions, where the BLUEs of the joint estimation, the OSP-based estimation, and the differential estimation are carried out for each one of the discussed time-based models. Some LS estimators without a proper whitening process are also presented for comparison. The acronyms of all estimators used in the simulations are summarized in Table 3. We also calculate the Cramér-Rao lower bound (CRLB) with an unknown $r_0$ based on the original model (24) [1, Chapter 3], since the TSE, SD-TOA, and SD-TDOA models all lose some information by ignoring some high-order terms. The root mean square error (RMSE) of the location estimate, which is defined as $\sqrt{E[(\hat{\mathbf{x}} - \mathbf{x})^2]}$ in general, is used as a performance measure in this paper. From the numerical results in Fig. 2, we can draw the following conclusions.

1. For each model, the corresponding BLUEs yield the same performance as expected.
2. Without a proper whitening, it can be observed that the performance of the LS estimators deteriorates. The D-LS-TSE-TOA, J-LS-SD-TOA, and J-LS-SD-TDOA clearly perform worse than their corresponding BLUEs.
3. The TSE model ignores $\mathcal{O}((\mathbf{x}_t - \hat{\mathbf{x}}_t^{(k-1)})^2)$ and accordingly suffers some information loss in modeling. However, the information loss can be reduced with a more accurate $\hat{\mathbf{x}}_t^{(k-1)}$. Therefore, with more iterations, the BLUEs for the TSE model approach the CRLB, which is in fact the essence of the ML property.
4. The SD-TOA model ignores $n_i^2, \forall i$ while the SD-TDOA model ignores $n_{i,j}^2, \forall i, i \neq j$. Ignoring these terms will cause an increasing information loss as the measurement noise gets larger.
5. Even though the BLUEs of the SD-TOA model outperform those of the SD-TDOA model in our simulation, we still cannot decide at this point which model is the best. This is because an optimal localization problem for the SD models should also include any dependence between the (nuisance) parameters, e.g., between $\mathbf{x}_t$ and $||\mathbf{x}_t||_2^2$, between $r_0$ and $r_0^2$ in $\boldsymbol{\theta}_1$, or between $\mathbf{x}_t$ and $r_j$ in $\boldsymbol{\theta}_2$, which explains the huge gap between the CRLB and the BLUEs for the SD models. By contrast, the TSE model obviously does not have this kind of issue. Nevertheless, including these dependencies is beyond the scope of this paper and we will not further consider this.
6. In practice, both the TSE and SD methods require iterations to obtain an accurate location estimate. However, note that, even after serveral iterations, the estimators based on the SD models still need to cope with the abovementioned dependency issue. Therefore, in real life, one often combines those two models, i.e., one uses the TSE model with the J-LS-SD-TDOA or the J-LS-SD-TOA as an initialization.
7. For the SD-TDOA model, ignoring the terms $n_{i,j}^2, \forall i, i \neq j$ implies that the information loss depends on the reference choice of the differencing process in (32). However, this is only because of the SD modeling thereafter, not because of the differencing process itself. Note that, for any other differencing process in this paper, the reference index is not important as long as the model is properly whitened.

**Table 3** Acronyms of the estimators used in the localization simulations

| Notations | Data models | Estimation methods |
| --- | --- | --- |
| J-BLUE-TSE-TOA, $k=1$ | White TSE model (25)[a], $M=1$ | Joint estimation (2) |
| OSP-BLUE-TSE-TOA, $k=1$ | " | OSP-based estimation (6) or (9) |
| D-BLUE-TSE-TOA[b], $k=1$ | " | Differential estimation (22) |
| D-LS-TSE-TOA[b], $k=1$ | " | LS estimator based on the unwhitened differential observations in (19) |
| J-LS-SD-TOA | Unwhitened SD-TOA model (29), $M=2$ | LS estimator with correlated model noise |
| J-BLUE-SD-TOA | Whitened SD-TOA model (31b), $M=2$ | Joint estimation (2) |
| OSP-BLUE-SD-TOA | " | OSP-based estimation (6) or (9) |
| D-BLUE-SD-TOA | " | Differential estimation (22) |
| J-LS-SD-TDOA | Unwhitened SD-TDOA model (34), $M=1$ | LS estimator with correlated model noise |
| J-BLUE-SD-TDOA | Whitened SD-TDOA model (36b), $M=1$ | Joint estimation (2) |
| OSP-BLUE-SD-TDOA | " | OSP-based estimation (6) or (9) |
| D-BLUE-SD-TDOA | " | Differential estimation (22) |
| J-LS-SD-RSS | Unwhitened SD-RSS model (40), $M=1$ | LS estimator with correlated model noise |
| J-BLUE-SD-RSS | Whitened SD-RSS model (43b), $M=1$ | Joint estimation (2) |
| OSP-BLUE-SD-RSS | " | OSP-based estimation (6) or (9) |
| D-LS-SD-RSS[c] | " | LS estimator based on the unwhitened differential observations in (19) |
| D-BLUE-SD-RSS[c] | " | Differential estimation (22) |

[a] The J-LS-SD-TDOA is used as an initial value (i.e., $k=0$), which is guaranteed to be near the global solution
[b] D-BLUE-TSE-TOA and D-LS-TSE-TOA can equivalently be considered to work with the TDOA measurements
[c] D-LS-SD-RSS and D-BLUE-SD-RSS can equivalently be considered to work with the DRSS measurements

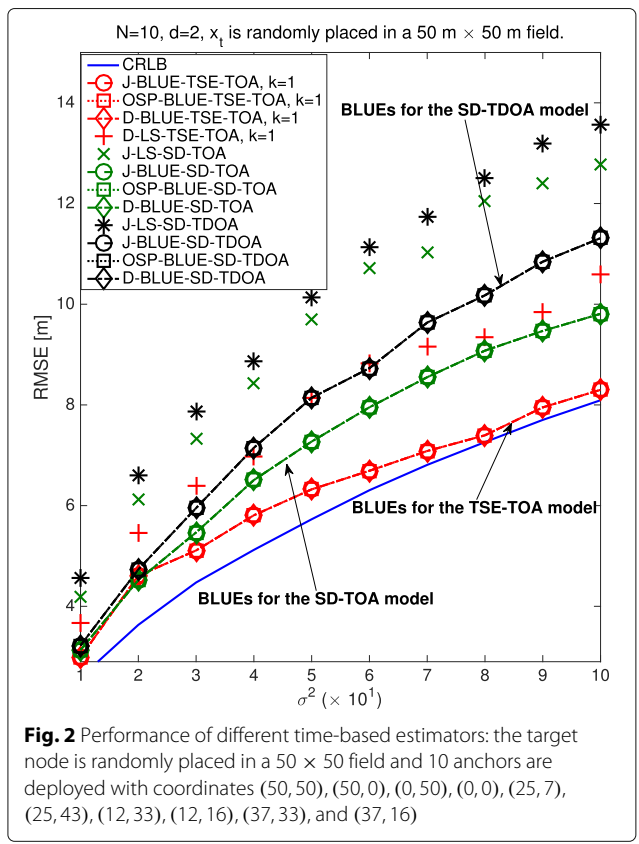

**Fig. 2** Performance of different time-based estimators: the target node is randomly placed in a 50 × 50 field and 10 anchors are deployed with coordinates (50, 50), (50, 0), (0, 50), (0, 0), (25, 7), (25, 43), (12, 33), (12, 16), (37, 33), and (37, 16)

### 3.2 Received signal strength based localization

Due to the simplicity of utilizing received signal strength (RSS) measurements, wireless networks with very constrained resources preferably rely on RSS-based localization [2]. Therefore, it gradually became very popular in recent years, and many efforts have already been put on this topic [55–58].

RSS-based localization mainly suffers from the complicated radio propagation channel. As before, assume that the target node is located at $\mathbf{x}_t$ and the $i$th anchor at $\mathbf{s}_i$. Based on a large-scale log-normal fading model [59], the RSS measurement can then be modeled as

$$P_i = P_0 - 10\gamma \log_{10}\left(\frac{||\mathbf{x}_t - \mathbf{s}_i||_2}{d_0}\right) + n_i, \ i = 1, 2, \cdots, N, \quad (37)$$

where $P_0$ is the received power at the reference distance $d_0$, $\gamma$ is the path-loss exponent (PLE), $n_i \sim \mathcal{N}(0, \sigma^2)$ is the shadowing effect, and $N$ is the number of anchor nodes. RSS-based localization is aimed at estimating the target location $\mathbf{x}_t$ from the RSS measurements. However, in some military or hostile scenarios, the transmit power might be unknown. Therefore, without loss of generality, we assume the reference distance $d_0$ to be 1 m and then the problem of the unknown transmit power can be equivalently converted into that of an unknown $P_0$. Note that (37) also has the non-linearity issue and, obviously, the iterative TSE model for RSS-based localization will be very similar to that developed for time-based localization. Therefore, to save space, we do not consider directly applying the TSE model but only focus on the SD method here.

To construct a linear data model, we rewrite (37) as

$$||\mathbf{x}_t - \mathbf{s}_i||_2^2 = \frac{P_0' n_i'}{P_i'}, \quad (38)$$

where $P_i' \triangleq 10^{\frac{P_i}{5\gamma}}, P_0' \triangleq 10^{\frac{P_0}{5\gamma}}$ and $n_i' \triangleq 10^{\frac{n_i}{5\gamma}}$. Interestingly though, we still need to apply the TSE to $n_i'$ here[2], such that (38) can further be approximated as

$$||\mathbf{x}_t||_2^2 - 2\mathbf{s}_i^T \mathbf{x}_t + ||\mathbf{s}_i||_2^2 = \frac{P_0'}{P_i'}\left[1 + \frac{ln(10)}{5\gamma} n_i\right]. \quad (39)$$

Then, a linear SD-RSS model for localization can be formulated from (39) as

$$\mathbf{h} = \mathbf{F}\boldsymbol{\phi} + \boldsymbol{\varsigma} \quad (40)$$

where

$$\mathbf{F} \triangleq \begin{bmatrix} \vdots & \vdots & \vdots \\ 2\mathbf{s}_i^T & -1 & 1/P_i' \\ \vdots & \vdots & \vdots \end{bmatrix}_{N \times (d+2)}, \quad (41a)$$

$$\boldsymbol{\phi} \triangleq \begin{bmatrix} \mathbf{x}_t \\ ||\mathbf{x}_t||_2^2 \\ P_0' \end{bmatrix}_{(d+2) \times 1}, \quad (41b)$$

$$\mathbf{h} \triangleq \begin{bmatrix} \vdots \\ ||\mathbf{s}_i||_2^2 \\ \vdots \end{bmatrix}_{N \times 1}, \quad (41c)$$

$$\boldsymbol{\varsigma} \triangleq \begin{bmatrix} \vdots \\ \frac{ln(10)P_0'}{5\gamma P_i'} n_i \\ \vdots \end{bmatrix}_{N \times 1}. \quad (41d)$$

This model was firstly presented in [57, eq. (18)] but in the absence of the shadowing effect. If we whiten the model (40) utilizing the covariance matrix of $\boldsymbol{\varsigma}$, i.e.,

$$\boldsymbol{\Sigma}_{\boldsymbol{\varsigma}} = \frac{[ln(10)]^2 P_0'^2 \sigma^2}{25\gamma^2} \mathbf{D}^{-2}, \quad (42)$$

where $\mathbf{D} = \text{diag}([P_1', \cdots, P_N']^T)$, we can obtain

$$\boldsymbol{\Sigma}_{\boldsymbol{\varsigma}}^{-1/2} \mathbf{h} = \boldsymbol{\Sigma}_{\boldsymbol{\varsigma}}^{-1/2} \mathbf{F}\boldsymbol{\phi} + \boldsymbol{\Sigma}_{\boldsymbol{\varsigma}}^{-1/2} \boldsymbol{\varsigma} \quad (43a)$$

$$\Rightarrow \mathbf{Dh} = \mathbf{DF}\boldsymbol{\phi} + \mathbf{D}\boldsymbol{\varsigma} \quad (43b)$$

where the covariance matrix of $\mathbf{D}\boldsymbol{\varsigma}$ becomes a scaled identity matrix, i.e., $\boldsymbol{\Sigma}_{\mathbf{D}\boldsymbol{\varsigma}} = \frac{ln(10)^2 P_0'^2 \sigma^2}{25\gamma^2} \mathbf{I}_N$. Note that this whitening step simply corresponds to an appropriate scaling of every entry of (40).

The whitened model (43b) is found to match our general model (1), since we notice that $\mathbf{DF}$ can be split into

$$\mathbf{DF} = \begin{bmatrix} \vdots & \vdots & \vdots \\ 2\mathbf{s}_i^T P_i' & -P_i' & 1 \\ \vdots & \vdots & \vdots \end{bmatrix} = \begin{bmatrix} \mathbf{F}' & \mathbf{1}_{N \times 1} \end{bmatrix}, \quad (44)$$

where $\mathbf{F}'$ contains the first $d+1$ columns of $\mathbf{DF}$. The relation between this model and the general model (1) is presented in Table 2. Note that we only consider a single nuisance parameter $P_0'$ in this model ($M = 1$). Although we could consider both $||\mathbf{x}_t||_2^2$ and $P_0'$ as nuisance parameters ($M = 2$), which would lead to the same performance after using the correct preprocessing steps, the reason why we take $M = 1$ here is to connect this model to the existing literature. For instance, after removing $P_0'$ using a single differencing step, the model for $\boldsymbol{\Gamma}_j \mathbf{Dh}$ is equal to the SD-DRSS model used in [57, eq. (22)]. However, without an appropriate whitening procedure, the LS estimators of the SD-RSS and SD-DRSS models yield a different performance, which is why they were treated and studied separately. Now, we realize that they actually are identical to each other as long as the model noise is properly whitened.

**Numerical results:** A simulation has also been conducted to verify our conclusions for this example. As before, the BLUEs of the joint estimation, OSP-based estimation, and the differential estimation for the SD-RSS model are evaluated and compared with some LS estimators without a proper whitening. Based on the original model in (37), the CRLB with an unknown $P_0$ is easy to calculate [1, Chapter 3]. From the numerical results in Fig. 3, the critical observation is that all the BLUEs here yield exactly the same performance as expected. Due to the colored model noise, the J-LS-SD-RSS and the D-LS-SD-RSS are relatively worse. Finally, denoting $R \triangleq ||\mathbf{x}_t||_2^2$, we again point out that neglecting the dependence between $R$ and $\mathbf{x}_t$ results in the gap between the CRLB and the estimators presented here.

### 3.3 Other examples
We believe that there are many other examples with linear nuisance parameters for our results. However, due to the limited space, we will only point out some of them. Besides the aforementioned localization examples, if anchors are separated into groups with different central clocks, multiple relative clock biases might exist in the TDOA measurements for localization, which can be removed by the OSP method [60, eq. (3)]. In cooperative localization, the multidimentional scaling (MDS) also uses the OSP-based method to eliminate the unknown terms [61, eq. (3)]. An acoustic source localization model, which also matches our general model (1), was presented in [62, eq. (6)].

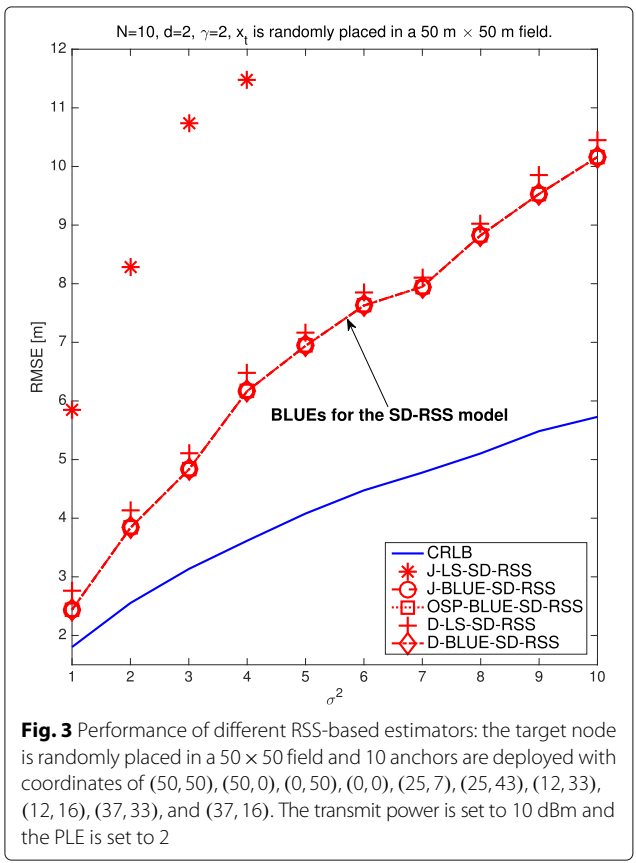

**Fig. 3** Performance of different RSS-based estimators: the target node is randomly placed in a 50 × 50 field and 10 anchors are deployed with coordinates of (50, 50), (50, 0), (0, 50), (0, 0), (25, 7), (25, 43), (12, 33), (12, 16), (37, 33), and (37, 16). The transmit power is set to 10 dBm and the PLE is set to 2

In [4, eq. (2)], the transmission times and clock offsets are the unknown nuisance parameters for the considered clock synchronization problem. The authors claim that those unknown parameters are systematically ML estimated before the synchronization. However, in fact, those nuisance parameters are equivalently removed by using respectively the observations $\mathbf{d}_{avg}$ in (13) or the OSP procedure. In hyperspectral imaging, OSP is also a very common procedure to extract the desired signals [19]. And when tracking mobile targets, frequency-difference-of-arrival measurements are often measured to cope with the Doppler effect [17, 18, 63, 64]. Furthermore, multiple-input-multiple-output (MIMO) receiver design might be affected by some nuisance parameters like I-Q imbalance and DC offset [5, eq. (7)]. In machine learning, a well-designed OSP is desired for dimensionality reduction [8, 9]. Extracting and working on the signal space is a strong need for signal separation [7] and underwater communication [6], which can be facilitated by OSP. At last, the famous differential global positioning system (DGPS) introduces a reference station on the ground and constructs a new differential observation set for positioning [65], where even the double differencing process is considered [66–68].

## 4 Conclusions

In this paper, we have introduced a general framework for estimation in the presence of unknown linear nuisance parameters. Three different kinds of methods to cope with the unknown nuisance parameters have been studied, i.e., the joint estimation, the OSP-based estimation, and the differential estimation. These approaches have been analyzed by investigating their corresponding BLUEs, where a new differential method has been introduced to cope with multiple nuisance parameters. We have discovered that, after a proper whitening procedure, all the BLUEs are equivalent to each other. From this interesting fact, one can draw some useful conclusions:

1. There only exists one unique BLUE for all these methods proposed to cope with unknown nuisance parameters.
2. Compared with the joint estimation, which directly utilizes all the original observations, none of the other two methods suffers any information loss.
3. For the differential approach, which requires selecting some references, the choice of the references is not important since there is no actual trace of the selected references in the corresponding BLUE.
4. In the differencing process, compared with the full differential observation set, any subset related to a single reference already preserves the full data information.

The presented analyses of the general model can be projected onto many practical applications, e.g., hyperspectral imaging, source localization and synchronization. Some localization examples have also been demonstrated, simulated and discussed to verify our conclusions.

### Endnotes

[1] For example, the noise **n** could also be uniform, Laplace, or student's *t*-distributed [69]

[2] We use $a^x = 1 + x\ln(a) + \cdots + \frac{(x\ln(a))^n}{n!} + \cdots, -\infty < x < \infty$ [70]. Note that the right hand side of (39) is an approximation, but it is regarded to be exact in this paper.

**References**
1. SM Kay, *Fundamentals of Statistical Signal Processing: Estimation Theory*. (Prentice-Hall, Inc., Upper Saddle River, NJ, USA, 1993)
2. N Patwari, JN Ash, S Kyperountas, AO Hero, RL Moses, NS Correal, Locating the nodes: cooperative localization in wireless sensor networks. IEEE Signal Proc. Mag. **22**(4), 54–69 (2005). doi:10.1109/MSP.2005.1458287
3. S Lopez, T Vladimirova, C Gonzalez, J Resano, D Mozos, A Plaza, The Promise of Reconfigurable Computing for Hyperspectral Imaging Onboard Systems: a Review and Trends. Proc. IEEE. **101**(3), 698–722 (2013). doi:10.1109/JPROC.2012.2231391

4. O Jean, AJ Weiss, Passive Localization and Synchronization Using Arbitrary Signals. IEEE Trans. Signal Process. **62**(8), 2143–2150 (2014). doi:10.1109/TSP.2014.2307281
5. CJ Hsu, R Cheng, WH Sheen, Joint least squares estimation of frequency, DC offset, I-Q imbalance, and channel in MIMO receivers. IEEE Trans. Veh. Technol. **58**(5), 2201–2213 (2009). doi:10.1109/TVT.2008.2005989
6. J He, MNS Swamy, MO Ahmad, Joint space-time parameter estimation for underwater communication channels with velocity vector sensor arrays. IEEE Trans. Wirel. Commun. **11**(11), 3869–3877 (2012). doi:10.1109/TWC.2012.092112.110875
7. MA Uusitalo, RJ Ilmoniemi, Signal-space projection method for separating MEG or EEG into components. Med. Biol. Eng. Comput. **35**(2), 135–140 (1997)
8. E Kokiopoulou, Y Saad, Orthogonal neighborhood preserving projections: a projection-based dimensionality reduction technique. IEEE Trans. Pattern Anal. Mach. Intell. **29**(12), 2143–2156 (2007). doi:10.1109/TPAMI.2007.1131
9. CX Ren, DQ Dai, in *Pattern Recognition, 2008. CCPR '08. Chinese Conference On*. 2d-onpp: Two dimensional extension of orthogonal neighborhood preserving projections for face recognition (IEEE, Beijing, 2008), pp. 1–6. doi:10.1109/CCPR.2008.48
10. S Bar, J Tabrikian, Bayesian estimation in the presence of deterministic nuisance parameters-part i: performance bounds. IEEE Trans. Signal Process. **63**(24), 6632–6646 (2015). doi:10.1109/TSP.2015.2468684
11. S Bar, J Tabrikian, Bayesian estimation in the presence of deterministic nuisance parameters-part ii: estimation methods. IEEE Trans. Signal Process. **63**(24), 6647–6658 (2015). doi:10.1109/TSP.2015.2468680
12. S Zhu, Z Ding, Joint synchronization and localization using TOAs: a linearization based WLS solution. IEEE J. Sel. Areas Commun. **28**(7), 1017–1025 (2010). doi:10.1109/JSAC.2010.100906
13. T-M Tu, C-H Chen, C-I Chang, A posteriori least squares orthogonal subspace projection approach to desired signature extraction and detection. IEEE Trans. Geosci. Remote Sens. **35**(1), 127–139 (1997). doi:10.1109/36.551941
14. YT Chan, KC Ho, A simple and efficient estimator for hyperbolic location. IEEE Trans. Signal Process. **42**(8), 1905–1915 (1994). doi:10.1109/78.301830
15. KC Ho, Bias Reduction for an Explicit Solution of Source Localization Using TDOA. IEEE Trans. Signal Process. **60**(5), 2101–2114 (2012). doi:10.1109/TSP.2012.2187283
16. Y Wang, G Leus, Reference-free time-based localization for an asynchronous target. EURASIP J. Adv. Signal Process. **2012**(1), 19 (2012). doi:10.1186/1687-6180-2012-19
17. KC Ho, X Lu, L Kovavisaruch, Source Localization Using TDOA and FDOA Measurements in the Presence of Receiver Location Errors: Analysis and Solution. IEEE Trans. Signal Process. **55**(2), 684–696 (2007). doi:10.1109/TSP.2006.885744
18. D Musicki, W Koch, in *Information Fusion, 2008 11th International Conference On*. Geolocation using TDOA and FDOA measurements (IEEE, Cologne, 2008), pp. 1–8
19. JC Harsanyi, C-I Chang, Hyperspectral image classification and dimensionality reduction: an orthogonal subspace projection approach. IEEE Trans. Geosci. Remote Sens. **32**(4), 779–785 (1994). doi:10.1109/36.298007
20. C-I Chang, Orthogonal subspace projection (OSP) revisited: a comprehensive study and analysis. IEEE Trans. Geosci. Remote Sens. **43**(3), 502–518 (2005). doi:10.1109/TGRS.2004.839543
21. M Song, CI Chang, A Theory of Recursive Orthogonal Subspace Projection for Hyperspectral Imaging. IEEE Trans. Geosci. Remote Sens. **53**(6), 3055–3072 (2015). doi:10.1109/TGRS.2014.2367816
22. Q Xu, Y Lei, J Cao, H Wei, in *Image and Signal Processing (CISP), 2014 7th International Congress On*. An improved algorithm based on reference selection for time difference of arrival location (IEEE, Dalian, 2014), pp. 953–957. doi:10.1109/CISP.2014.7003916
23. Y Wang, F Zheng, M Wiemeler, W Xiong, T Kaiser, in *Vehicular Technology Conference (VTC Fall), 2013 IEEE 78th*. Reference Selection for Hybrid TOA/RSS Linear Least Squares Localization (IEEE, Las Vegas, 2013), pp. 1–5. doi:10.1109/VTCFall.2013.6692388
24. I Guvenc, S Gezici, F Watanabe, H Inamura, in *2008 IEEE Wireless Communications and Networking Conference*. Enhancements to Linear Least Squares Localization Through Reference Selection and ML Estimation (IEEE, Las Vegas, 2008), pp. 284–289. doi:10.1109/WCNC.2008.55
25. HC So, YT Chan, FKW Chan, Closed-Form Formulae for Time-Difference-of-Arrival Estimation. IEEE Trans. Signal Process. **56**(6), 2614–2620 (2008). doi:10.1109/TSP.2007.914342
26. SCK Herath, PN Pathirana, Robust Localization With Minimum Number of TDoA Measurements. IEEE Signal. Proc. Let. **20**(10), 949–951 (2013). doi:10.1109/LSP.2013.2274273
27. Y Huang, J Benesty, GW Elko, RM Mersereati, Real-time passive source localization: a practical linear-correction least-squares approach. IEEE T. Speech. Audi. P. **9**(8), 943–956 (2001). doi:10.1109/89.966097
28. RO Schmidt, A New Approach to Geometry of Range Difference Location. IEEE Trans. Aerosp. Electron. Syst. **AES-8**(6), 821–835 (1972). doi:10.1109/TAES.1972.309614
29. R Schmidt, Least squares range difference location. IEEE Trans. Aerosp. Electron. Syst. **32**(1), 234–242 (1996). doi:10.1109/7.481265
30. S Venkatesh, RM Buehrer, in *Proceedings of the 5th International Conference on Information Processing in Sensor Networks*. IPSN '06. A Linear Programming Approach to NLOS Error Mitigation in Sensor Networks (ACM, New York, NY, USA, 2006), pp. 301–308. doi:10.1145/1127777.1127823. http://doi.acm.org/10.1145/1127777.1127823
31. A-J van der Veen, EF Deprettere, AL Swindlehurst, Subspace-based signal analysis using singular value decomposition. Proc. IEEE. **81**(9), 1277–1308 (1993). doi:10.1109/5.237536
32. LL Scharf, ML McCloud, Blind adaptation of zero forcing projections and oblique pseudo-inverses for subspace detection and estimation when interference dominates noise. IEEE Trans. Signal Process. **50**(12), 2938–2946 (2002). doi:10.1109/TSP.2002.805245
33. AH Sayed, A Tarighat, N Khajehnouri, Network-based wireless location: challenges faced in developing techniques for accurate wireless location information. IEEE Signal Proc. Mag. **22**(4), 24–40 (2005). doi:10.1109/MSP.2005.1458275
34. WH Foy, Position-Location Solutions by Taylor-Series Estimation. IEEE Trans. Aerosp. Electron. Syst. **AES-12**(2), 187–194 (1976). doi:10.1109/TAES.1976.308294
35. CT Kelley, Iterative Methods for Optimization. Front. Appl. Math. Soc. Ind. Appl. Math (1999). https://books.google.nl/books?id=Bq6VcmzOe1IC
36. A Beck, P Stoica, J Li, Exact and Approximate Solutions of Source Localization Problems. IEEE Trans. Signal Process. **56**(5), 1770–1778 (2008). doi:10.1109/TSP.2007.909342
37. DB Haddad, WA Martins, MdVM da Costa, LWP Biscainho, LO Nunes, B Lee, Robust Acoustic Self-Localization of Mobile Devices. in *IEEE Transactions on Mobile Computing*. **15**(4), 982–995 (2016). doi:10.1109/TMC.2015.2439278
38. KW Cheung, HC So, WK Ma, YT Chan, Least squares algorithms for time-of-arrival-based mobile location. IEEE Trans. Signal Process. **52**(4), 1121–1130 (2004). doi:10.1109/TSP.2004.823465
39. JC Chen, RE Hudson, K Yao, Maximum-likelihood source localization and unknown sensor location estimation for wideband signals in the near-field. IEEE Trans. Signal Process. **50**(8), 1843–1854 (2002). doi:10.1109/TSP.2002.800420
40. C-H Park, S Lee, J-H Chang, Robust closed-form time-of-arrival source localization based on alpha-trimmed mean and HodgesCLehmann estimator under NLOS environments. Signal Process. **111**, 113–123 (2015). doi:10.1016/j.sigpro.2014.12.020
41. M Sun, KC Ho, Successive and Asymptotically Efficient Localization of Sensor Nodes in Closed-Form. IEEE Trans. Signal Process. **57**(11), 4522–4537 (2009). doi:10.1109/TSP.2009.2025821
42. ND Gaubitch, WB Kleijn, R Heusdens, in *2013 IEEE International Conference on Acoustics, Speech and Signal Processing*. Auto-localization in ad-hoc microphone arrays (IEEE, Vancouver, 2013), pp. 106–110. doi:10.1109/ICASSP.2013.6637618
43. L Wang, TK Hon, JD Reiss, A Cavallaro, Self-Localization of Ad-Hoc Arrays Using Time Difference of Arrivals. IEEE Trans. Signal Process. **64**(4), 1018–1033 (2016). doi:10.1109/TSP.2015.2498130
44. K Liu, X Liu, X Li, Guoguo: Enabling Fine-Grained Smartphone Localization via Acoustic Anchors. IEEE Trans. Mob. Comput. **15**(5), 1144–1156 (2016). doi:10.1109/TMC.2015.2451628
45. JJ Caffery, in *Vehicular Technology Conference, 2000. IEEE-VTS Fall VTC 2000. 52nd*. A new approach to the geometry of TOA location, vol. 4 (IEEE, Boston, 2000), pp. 1943–19494. doi:10.1109/VETECF.2000.886153

46. Y Wang, G Leus, X Ma, in *2011 IEEE International Conference on Acoustics, Speech and Signal Processing (ICASSP)*. Time-based localization for asynchronous wireless sensor networks (IEEE, Prague, 2011), pp. 3284–3287. doi:10.1109/ICASSP.2011.5946723
47. P Stoica, J Li, Lecture Notes - Source Localization from Range-Difference Measurements. IEEE Signal Process. Mag. **23**(6), 63–66 (2006). doi:10.1109/SP-M.2006.248717
48. Y Liu, F Guo, L Yang, W Jiang, An Improved Algebraic Solution for TDOA Localization With Sensor Position Errors. IEEE Commun. Lett. **19**(12), 2218–2221 (2015). doi:10.1109/LCOMM.2015.2486769
49. J Liu, Z Wang, JH Cui, S Zhou, B Yang, A Joint Time Synchronization and Localization Design for Mobile Underwater Sensor Networks. IEEE Trans. Mob. Comput. **15**(3), 530–543 (2016). doi:10.1109/TMC.2015.2410777
50. B Huang, L Xie, Z Yang, TDOA-Based Source Localization with Distance-Dependent Noises. IEEE Trans. Wirel. Commun. **14**(1), 468–480 (2015). doi:10.1109/TWC.2014.2351798
51. H Yang, J Chun, D Chae, Hyperbolic Localization in MIMO Radar Systems. IEEE Antennas Wirel. Propag. Lett. **14**, 618–621 (2015). doi:10.1109/LAWP.2014.2374603
52. J Smith, J Abel, The spherical interpolation method of source localization. IEEE J. Oceanic. Eng. **12**(1), 246–252 (1987). doi:10.1109/JOE.1987.1145217
53. B Friedlander, A passive localization algorithm and its accuracy analysis. IEEE J. Oceanic. Eng. **12**(1), 234–245 (1987). doi:10.1109/JOE.1987.1145216
54. J Smith, J Abel, Closed-form least-squares source location estimation from range-difference measurements. IEEE Trans. Acoust. Speech Signal Process. **35**(12), 1661–1669 (1987). doi:10.1109/TASSP.1987.1165089
55. X Li, RSS-Based Location Estimation with Unknown Pathloss Model. IEEE Trans. Wirel. Commun. **5**(12), 3626–3633 (2006). doi:10.1109/TWC.2006.256985
56. HC So, L Lin, Linear Least Squares Approach for Accurate Received Signal Strength Based Source Localization. IEEE Trans. Signal Process. **59**(8), 4035–4040 (2011). doi:10.1109/TSP.2011.2152400
57. RM Vaghefi, MR Gholami, EG Strom, in *Acoustics, Speech and Signal Processing (ICASSP), 2011 IEEE International Conference On*. RSS-based sensor localization with unknown transmit power (IEEE, Prague, 2011), pp. 2480–2483. doi:10.1109/ICASSP.2011.5946987
58. MR Gholami, RM Vaghefi, EG Strom, RSS-Based Sensor Localization in the Presence of Unknown Channel Parameters. IEEE Trans. Signal Process. **61**(15), 3752–3759 (2013). doi:10.1109/TSP.2013.2260330
59. T Rappaport, *Wireless Communications: Principles and Practice*, 2nd edn. (Prentice Hall PTR, Upper Saddle River, NJ, USA, 2001)
60. Y Wang, KC Ho, TDOA Source Localization in the Presence of Synchronization Clock Bias and Sensor Position Errors. IEEE Trans. Signal Process. **61**(18), 4532–4544 (2013). doi:10.1109/TSP.2013.2271750
61. S Kumar, R Kumar, K Rajawat, Cooperative localization of mobile networks via velocity-assisted multidimensional scaling. IEEE Trans. Signal Process. **64**(7), 1744–1758 (2016). doi:10.1109/TSP.2015.2507548
62. D Li, YH Hu, in *Parallel Processing Workshops, 2004. ICPP 2004 Workshops. Proceedings. 2004 International Conference On*. Least square solutions of energy based acoustic source localization problems (IEEE, Montreal, 2004), pp. 443–446. doi:10.1109/ICPPW.2004.1328053
63. HW Wei, R Peng, Q Wan, ZX Chen, SF Ye, Multidimensional Scaling Analysis for Passive Moving Target Localization With TDOA and FDOA Measurements. IEEE Trans. Signal Process. **58**(3), 1677–1688 (2010). doi:10.1109/TSP.2009.2037666
64. KC Ho, W Xu, An accurate algebraic solution for moving source location using TDOA and FDOA measurements. IEEE Trans. Signal Process. **52**(9), 2453–2463 (2004). doi:10.1109/TSP.2004.831921
65. BW Parkinson, JJ Spilker, Global Positioning System: Theory and Applications. Progress in astronautics and aeronautics. Am. Inst. Aeronaut. Astronaut. **v. 1** (1996). https://books.google.nl/books?id=lvI1a5J_4ewC
66. RO Nielsen, Relationship between dilution of precision for point positioning and for relative positioning with GPS. IEEE Trans. Aerosp. Electron. Syst. **33**(1), 333–338 (1997). doi:10.1109/7.570809
67. PJG Teunissen, A proof of Nielsen's conjecture on the GPS dilution of precision. IEEE Trans. Aerosp. Electron. Syst. **34**(2), 693–695 (1998). doi:10.1109/7.670364
68. C Park, I Kim, Comments on "relationships between dilution of precision for point positioning and for relative positioning with GPS". IEEE Trans. Aerosp. Electron. Syst. **36**(1), 315–316 (2000). doi:10.1109/7.826336
69. M Abramowitz, IA Stegun, et al, Handbook of mathematical functions. Appl. Math. Ser. **55**, 62 (1966)
70. M Abramowitz, *[Handbook of Mathematical Functions, With Formulas, Graphs, and Mathematical Tables]*. (Dover Publications, Incorporated, Mineola, 1974)